\documentclass[reqno,centertags, 12pt]{amsart}
\numberwithin{equation}{section}

\newcommand{\vp}{\varphi}
\newcommand{\T}{\partial\mathbb{D}}

\newcommand{\ds}{\displaystyle}
\newcommand{\ol}{\overline}

\newcommand{\esssupp}{\text{\rm{ess supp}}}

\newcommand{\be}{\begin{equation}}
\newcommand{\ee}{\end{equation}}
\newcommand{\ba}{\begin{array}}
\newcommand{\ea}{\end{array}}
\newcommand{\C}{\mathbb{C}}

\newcommand{\bpm}{\begin{pmatrix}}
\newcommand{\epm}{\end{pmatrix}}
\newcommand{\ainfty}{A_{\infty}}
\newcommand{\w}{\bpm  w_1 \\  w_2 \epm}
\newcommand{\sumn}{\sum_{n=0}^{\infty}}

\newcommand{\mT}{\mathcal{T}}
\newcommand{\mS}{\mathcal{S}}

\newcommand{\ye}{\left( \ds \frac y 2 \right)}
 
\newcommand{\yy}{\bpm 1 \\ 1 \epm}
\newcommand{\ninn}{n \in \mathbb{N}} 
\newcommand{\tr}{\text{Tr}}

\newtheorem{lemma}{Lemma}[section]
\newtheorem{theorem}{Theorem}[section]
\newtheorem{corollary}{Corollary}[section]

\begin{document}

\title[Asymptotics of OP's and point perturbation on the unit circle]
{Asymptotics of orthogonal polynomials and point perturbation on the unit circle}
\author[M.-W. L. Wong]{Manwah Lilian Wong}
\thanks{$^*$ Mathematics 253-37, California Institute of Technology, Pasadena, CA 91125.
E-mail: math.wmw@gmail.com}
\date{August 2, 2009}

\maketitle

\begin{abstract}
In the first five sections, we deal with the class of probability measures with asymptotically periodic Verblunsky coefficients of $p$-type bounded variation. The goal is to investigate the perturbation of the Verblunsky coefficients when we add a pure point to a gap of the essential spectrum.

For the asymptotically constant case, we give an asymptotic formula for the orthonormal polynomials in the gap, prove that the perturbation term converges and show the limit explicitly. Furthermore, we prove that the perturbation is of bounded variation. Then we generalize the method to the asymptotically periodic case and prove similar results.

In the last two sections, we show that the bounded variation condition can be removed if a certain symmetry condition is satisfied. Finally, we consider the special case when the Verblunsky coefficients are real with the rate of convergence being $c_n$. We prove that the rate of convergence of the perturbation is in fact $O(c_n)$. In particular, the special case $c_n=1/n$ will serve as a counterexample to the possibility that the convergence of the perturbed Verblunsky coefficients should be exponentially fast when a point is added to a gap.

\end{abstract}

\section{Introduction}
\subsection{Background}
Suppose $d\mu$ is a probability measure on the unit circle $\T=\{z \in \C: |z|=1\}$. We define an inner product and a norm on $L^2(\T, d\mu)$ respectively as follows:
\begin{eqnarray}
\left\langle f ,g \right\rangle & = & \ds \int_{\T} \ol{f(e^{i \theta})} g(e^{i \theta}) d\mu(\theta) \\
\|f\|_{d\mu} & = & \left( \ds \int_{\T} |f(e^{i \theta})|^2 d\mu(\theta) \right)^{1/2}
\end{eqnarray}

Using the inner product defined above, we can orthogonalize $1, z, z^2, \dots$ to obtain the family of monic orthogonal polynomials associated with the measure $d\mu$, namely, $(\Phi_n(z, d\mu))_{\ninn}$. We denote the normalized family as $(\vp_n(z, d\mu))_{\ninn}$.

Closely related to $\Phi_n(z)$ is the family of reversed polynomials, defined as $\Phi_n^*(z)=z^n \ol{\Phi_n(1/\ol{z})}$. They obey the well-known Szeg\H o recursion relation
\be
\Phi_{n+1}(z) = z \Phi_n (z)- \ol{\alpha_n} \Phi_n^*(z)
\label{eq01}
\ee and $\alpha_n$ is known as the $n$-th Verblunsky coefficient. The Szeg\H o recursion relations for the normalized families are
\begin{align}
\vp_{n+1}(z) & = (1-|\alpha_n|^2)^{-1/2}(z \vp_n(z)-\ol{\alpha_n} \vp_n^*(z)) \label{normrec1} \\
\vp_{n+1}^*(z) & = (1-|\alpha_n|^2)^{-1/2}(\vp_n^*(z) - \alpha_n z \vp_n(z)) \label{normrec2}
\end{align}

These recursion relations will be useful later in this paper. For more on orthogonal polynomials on the unit circle, the reader may refer to \cite{geronimus1, onefoot, simon1, simon2, szego}.

\subsection{The point mass problem}
We add a point mass $\zeta = e^{i\omega} \in \T$ with weight $0<\gamma<1$ to $d\mu$ in the following manner:
\be
d\nu = (1-\gamma) d\mu + \gamma \delta_\omega
\label{dnudef}
\ee  Our goal is to investigate $\alpha_n(d\nu)$. 

\textbf{Remark about notation:} From now on, any object without the label $(d\nu)$ is considered to be associated with the original measure $d\mu$, unless otherwise stated.

Point mass perturbation has a long history (see the Introduction of \cite{wong1}). One of the classic results is the following theorem by Geronimus \cite{geronimus1, geronimus2}:
\begin{theorem}(Geronimus) Suppose the probability measure $d\nu$ is defined as in (\ref{dnudef}). Then
\be \Phi_{n}(z, d\nu) = \Phi_n(z) - \ds \frac{\vp_n(z) K_{n-1}(z, \zeta)}{(1-\gamma) \gamma^{-1} + K_{n-1}(\zeta, \zeta)}
 \label{geronimus}
\ee where \be K_{n}(z,\zeta)  = \ds \sum_{j=0}^{n} \ol{\vp_j(\zeta)} \vp_j(z) \label{kndef} \ee and all objects without the label $(d\nu)$ are associated with the measure $d\mu$. 
\end{theorem}

Since $\Phi_n(0) = -\ol{\alpha_{n-1}}$, by putting $z=0$ into (\ref{geronimus}) one gets a formula relating the Verblunsky coefficients of $d\mu$ and $d\nu$.

Formula (\ref{geronimus}) was rediscovered by Nevai \cite{nevai1} for OPRL and by Cachafeiro--Marcell\'an \cite{cm1, cm2, cm3} for OPUC. For general measures on $\mathbb{C}$, the formula is from Cachafeiro--Marcell\'an \cite{cm4, cm5}. Using a totally different approach, Simon \cite{simon2} found the following formula for OPUC:
\be
\alpha_{n}(d\nu)=\alpha_{n}-q_{n}^{-1} \gamma \ol{\vp_{n+1}(\zeta)} \left( \ds \sum_{j=0}^{n} \alpha_{j-1} \frac{\|\Phi_{n+1}\|}{\|\Phi_j\|}\vp_j(\zeta) \right)
\label{simonformula}
\ee where $q_{n}  =  (1-\gamma) + \gamma K_{n}(\zeta,\zeta) ; \alpha_{-1}  = -1$.

In \cite{wong1, wong2}, we applied the Christoffel--Darboux formula to (\ref{simonformula}) and proved the following formula for $\alpha_n(d\nu)$:
\be
\alpha_n(d\nu) = \alpha_n(d\mu) + \Delta_n(\zeta)
\label{addpoint}
\ee where
\be
\Delta_n(\zeta) = \ds \frac{(1-|\alpha_n|^2)^{1/2} \ol{\vp_{n+1}(\zeta)}\vp_n^*(\zeta)}{(1-\gamma)\gamma^{-1} + K_n(\zeta,\zeta)}; \quad K_n(\zeta,\zeta) = \ds \sum_{j=0}^{n} |\vp_j(\zeta)|^2
\label{deltandef}
\ee This prompted us to study the asymptotic behavior of $\vp_n(z)$ on $\T$ in order to understand the asymptotics of (\ref{deltandef}).

In \cite{wong1}, we considered the class of probability measures with $\ell^2$ Verblunsky coefficients of bounded variation, i.e., 
\be
\ds \sumn |\alpha_n|^2<\infty \quad \text{ and } \quad \ds \sumn |\alpha_n - \alpha_{n+1} | < \infty .
\label{condell2} \ee

In this paper, we consider the class of measures with asymptotically periodic Verblunsky coefficients of $p$-type bounded variation (this term was first introduced in \cite{ps3}), i.e., given a periodic sequence $\beta_n$ of period $p$,
\be
\ds \lim_{n \to \infty} (\alpha_n - \beta_n) = 0 \text{ and }
\ds \sum_{n=0}^{\infty} |\alpha_{n+p} - \alpha_n| < \infty .
\label{eqn31a}
\ee

First, we handle the special case $p=1$; then we generalize the method to any $p$. It is well-known that any measure satisfying (\ref{eqn31a}) has the same essential spectrum as $d\mu_\beta$ (the measure associated with $(\beta_n)_{\ninn}$) which is supported on a finite number of bands. The reader may refer to Chapter 11 of \cite{simon2} for a detailed discussion of periodic Verblunsky coefficients.

\subsection{Gaps and Periodicity}
Before we move on to stating the results, it would be helpful to have a brief discussion about gaps and periodicity.

By an application of Weyl's Theorem to the CMV matrix (see Theorem 4.3.5 of \cite{simon1}), $\alpha_n \to L$ implies that $d\mu$ has the same essential spectrum as the measure $d\mu_0$ with Verblunsky coefficients $\alpha_n(d\mu_0) \equiv L$ (the measure $d\mu_0$ is known to be associated with the Geronimus polynomials). Besides, it is known that $d\mu_0$ is supported on the arc 
\be \Gamma_{L} = [\theta_{|L|}, 2\pi - \theta_{|L|}] \label{gammaldef}
\ee
where $\theta_{|L|} = 2 \arcsin(|L|)$, and $d\mu_0$ admits at most one single pure point in $[-\theta_{|L|}, \theta_{|L|}]$. In other words, there is a gap $\mathbb{G}_L$ in the spectrum, with at most one pure point inside. The reader may refer to Example 1.6.12 of \cite{simon1} for a detailed discussion.

Note that $\alpha_n \equiv L$ can be seen as a periodic sequence of period 1, in fact, there is a more general result concerning gaps in the spectrum for measures with periodic Verblunsky coefficients. The precise statement reads as follows (see Theorem 11.1.2 of \cite{simon2}): let $(\beta_n)_n$ be a periodic family of Verblunsky coefficients of period $p$, i.e., $\beta_n = \beta_{n+p}$ for all $n$. Let $d\mu_\beta$ be the associated measure.  Then $\{ e^{i \theta}: |\tr(T_p(e^{i \theta})| \leq 2\}$ is a closed set which is the union of $p$ closed intervals $B_1, \dots, B_p$ (which can only overlap at the endpoints). Let
\be
B = \ds \cup_{j=1}^{p} B_j
\label{bdef}
\ee
Moreover, $B$ is the essential support of the a.c. spectrum. In each disjoint open interval on $\T \backslash B$, $d\mu$ has either no support or a single pure point.

As a result, in both cases that we consider, there are gaps in the spectrum and when $z \in \T$ is in one of those open gaps, we have $|\tr \,T_p(z) |>2$.

The reader may refer to Chapter 11 of \cite{simon2} for a detailed discussion of periodic Verblunsky coefficients.

\section{Results}

First, we present a new method for computing the asymptotics of $\vp_n(z)$ in the gap of the spectrum when the family $(\alpha_n)_{\ninn}$ is asymptotically constant and of bounded variation (see formulae (\ref{vpn}) and (\ref{vpn*})). Applying that to the point mass formula (\ref{deltandef}), we prove the following result:

\begin{theorem} Let $(\alpha_n)_{\ninn}$ be the Verblunsky coefficients of the probability measure $d\mu$ on $\T$ such that
\begin{eqnarray}
\alpha_n \to L \in \mathbb{D} \backslash \{0\} \label{condition1}\\
\ds \sum_{j=0}^{\infty} |\alpha_{j+1} - \alpha_j | < \infty \label{condition2}
\end{eqnarray} 
Let $\mathbb{G}_L$ be the gap of the essential spectrum (not including the endpoints). We add a pure point $\zeta=e^{i \theta} \in \mathbb{G}_L$ to $d\mu$ to form $d\nu$ as in (\ref{dnudef}). Then one of the following is true: \\
\begin{enumerate}
\item If $\mu(\zeta) >0$, then the three sequences $(|\vp_n(\zeta)|)_{\ninn}$, $(\Delta_n(\zeta))_{\ninn}$ and $(\alpha_n(d\nu)-\alpha_n(d\mu))_{\ninn}$ tend to zero exponentially fast.
\item If $\mu(\zeta) = 0$, then \\
\begin{enumerate}
\item $\lim_{n \to \infty} \Delta_n(\zeta)$ exists, and \be
\Delta_\infty(\zeta) \equiv \ds \lim_{n \to \infty} \Delta_n(\zeta) = \ds \overline{h(\zeta)^{1/2}} \left[ \ds \frac{(\zeta-1) - h(\zeta)^{1/2}}{2\ol{L}} \right]
\label{limit}
\ee where
\be
h(\zeta) = (\zeta-1)^2 + 4\zeta |L|^2
\ee and we choose the branch of logarithm such that $(1)^{1/2} = 1$.
\item Furthermore, $|\Delta_\infty(\zeta) + L| = |L|$ and
\be 
\ds \lim_{n \to\ \infty} \alpha_n(d\nu) = L e^{i \omega}
\ee where
\begin{eqnarray}
\cos \omega & =  
 & \ds \frac{2  \sin^2 \left( \frac \theta 2 \right)-|L|^2}{|L|^2} \label{phase1} \\
\sin \omega &= 
 &\frac{ 2 \sin\left( \frac \theta 2 \right) \sqrt{|L|^2 -\sin^2\left(\frac{\theta}{2} \right) } }{|L|^2}  \label{phase2}
\end{eqnarray} 

\item  $(\Delta_n(\zeta))_{\ninn}$ is of bounded variation, i.e.,
\be
\ds \sum_{n=0}^{\infty} |\Delta_{n+1}(\zeta) - \Delta_n(\zeta) |< \infty
\ee \\

\end{enumerate}
\end{enumerate}
\label{theorem1}
\end{theorem}

\noindent \textit{Three remarks about Theorem \ref{theorem1}: \\
(i) Since $\alpha_n \to L \not = 0$, this measure has the same essential spectrum as the measure $d\mu_0$ with Verblunsky coefficients $\alpha_n(d\mu_0) \equiv L$, which is supported on the arc $\Gamma_{|L|}$ as defined in (\ref{gammaadef}). \\
(ii) Case (1) is a special case of Corollary 24.3 of \cite{simon3}, where Simon proved that varying the weight of an isolated pure point in the gap will result in exponentially small perturbation to $\alpha_n(d\mu)$. \\
(iii) By (2c), adding a pure point to the gap will preserve the bounded variation property of $(\alpha_n)_{\ninn}$. Hence, we can add a finite number of points inductively and generalize the result to finitely many pure points in the gap.}

\vspace{1cm}

Next, we will generalize the technique developed in the proof of Theorem \ref{theorem1} and prove the following result about measures with asymptotically periodic Verblunsky coefficients:

\begin{theorem}
Let $(\beta_n)_{\ninn}$ be a periodic family of Verblunsky coefficients of period $p$, i.e., $\beta_n = \beta_{n+p}$ for all $n$, and let $d\mu_\beta$ be the measure associated with it. Let $\Gamma_\beta$ be the union of open arcs which are the interiors of the bands that form $\esssupp (d\mu_\beta)$. Suppose the measure $d\mu$ has Verblunsky coefficients $(\alpha_n)_{\ninn}$ that are asymptotically $p$-periodic of bounded variation, i.e.,
\begin{eqnarray}
\ds \lim_{n \to \infty} (\alpha_n - \beta_n) = 0 , \label{condap1} \\
\ds \sum_{n=0}^{\infty} |\alpha_{n+p} - \alpha_n| < \infty . \label{condap2}
\end{eqnarray}

Now we add a pure point $\zeta \in \T \setminus \overline \Gamma_\beta$ to $d\mu$ as in (\ref{addpoint}). Then one of the following is true:
\begin{enumerate}
\item  $\mu(\zeta)>0$, then for each fixed $0 \leq j < p$, $\lim_{k \to \infty} \Delta_{kp+j}(\zeta) = 0$ exponentially fast. 
\item $\mu(\zeta) = 0$, then for each fixed $0 \leq j < p$, $\lim_{k \to \infty} \Delta_{kp+j} (\zeta)$ exists and
\be
\ds \sum_{k=0}^{\infty} |\Delta_{(k+1)p+j}(\zeta) - \Delta_{kp+j}(\zeta) | < \infty .
\ee
\end{enumerate}
\label{theorem2}
\end{theorem}

\noindent \textit{Remark about Theorem \ref{theorem2}: it is worth noting that if one adds a pure point $\zeta$ as in \eqref{addpoint} to the support $\Gamma_\beta$, then $\lim_{n \to \infty} \Delta_n(\zeta) =0$. This result was proven by Peherstorfer--Steinbauer (see Theorem 3 of \cite{ps3}).}

\medskip
Then we will prove the following result where $(\alpha_n)_{\ninn}$ is not necessarily of bounded variation:

\begin{theorem} Let $\zeta \in \T$ and $\mu(\zeta)=0$. Suppose $ \lim_{n \to \infty} \zeta^n \alpha_n = L$. Then
\be \lim_{n \to \infty} \zeta^n \Delta_n(\zeta) = -2L .
\ee As a result,
\be
\lim_{n \to \infty} \zeta^n \alpha_n(d\nu) = - \lim_{n \to \infty} \zeta^n \alpha_n(d\mu).
\ee
\label{theorem3}
\end{theorem}

Finally, we use Theorem \ref{theorem3} to prove Corollary \ref{corollary1} below to illustrate the non-exponential rate of convergence of $\Delta_n(\zeta)$ towards its limit. One might have guessed that the convergence should be exponentially fast, but we will show that this is not the case!

\begin{corollary} Let $\alpha_n = L + {c_n}$, where $L<0$, $c_n \in \mathbb{R}$ and $c_n \to 0$. Then
\be
\Delta_n(1) = -2L -  2 {c_n} + o\left( {c_n} \right) .
\ee
\label{corollary1}
\end{corollary}

In particular, when $c_n = 1/n$, we have the rate of convergence being $O(1/n)$ which is not exponential.

The reader may also refer to \cite{wong3} in which Wong demonstrated that point perturbation of a certain class of measures on the real line would result in non-exponential perturbation of the recurrence coefficients.

\vspace{1cm}

There are many papers about measures supported on an interval/arc, and about the perturbation of orthogonal polynomials with periodic recursion coefficients. For example, the reader may refer to \cite{bello, denisov,  peherstorfereast, pu, barrios, alfaro, dks}.

Bello-L\'opez \cite{bello} extended the well-known work of Rakhmanov \cite{rakhmanov1, rakhmanov2, rakhmanov3} and proved the following: let $0<a<1$ and $\theta_a= 2 \arcsin(a)$. If $d\mu$ is supported on the arc 
\be \Gamma_a=\{\zeta \in \T: |\arg(\zeta)|>\theta_a \}
\label{gammaadef}
\ee  such that the absolutely continuous part $w(\theta)>0$ on $\Gamma_a$, then $\lim_{n\to \infty} |\alpha_n|=a$. Bellos-L\'opez's result is restricted to measures that are absolutely continuous on the arc, and it was later extended to measures with infinitely many mass points outside the a.c.\ part of the support (see for example, \cite{barrios} and Theorem 13.4.4 of \cite{simon2}). However, unlike Theorem \ref{theorem1}, these results do not tell us whether $\Delta_n(\zeta)$ approaches a single point.

In \cite{peherstorfereast}, Peherstorfer--Steinbauer considered the situation where $d\mu$ is an absolutely continuous measure on $supp(d\mu)= \Gamma_a$ with the a.c. part $w(\theta)$ satisfying the Szeg\H o condition on $\Gamma_a$, i.e., 
\be
\ds \int_{\Gamma_a} \log w(\theta) \ds \frac{\sin(\frac{\theta}{2})}{\sqrt{\cos^2 (\frac{\theta_{|a|}}{2} ) - \cos^2 (\frac{\theta}{2}})}d\theta > -\infty .
\label{szegocondition}
\ee They proved that if we add a finite number of pure points to the gap to form the measure to $d\tau$, then $\lim_{n \to \infty} \alpha_n(d\tau)$ exists and the limit has norm $|a|$. In the Appendix, we are going to work out an example that demonstrates the existence of a large class of measures with Verblunsky coefficients $\alpha_n \to L$ of bounded variation that fail the Szeg\H o condition (\ref{szegocondition}).


Given such a result for orthogonal polynomials on the unit circle, one would expect a similar result for the real line. In \cite{pu}, Peherstorfer--Yuditskii gave the following result: for any Jacobi matrix $J$ whose spectrum is a finite gap set with the a.c.\ part of the spectral measure satisfying the Szeg\H o condition, then there is a unique Jacobi matrix $J_\infty$ in the isospectral torus such that the orthogonal polynomials of $J$ and $J_\infty$ have the same asymptotics away from the spectrum as $n \to\infty$. In particular, this implies that the Jacobi parameters of $J$ converge to the parameters of $J_\infty$ as $n \to \infty$.


\section{Tools}
For the convenience of the reader, a brief discussion of two major tools used in the proofs will be presented here.

\subsection{The Ces\`aro--Stolz Theorem}
One of the very important tools for the computation of the limit $\lim_{n \to \infty}\Delta_n(\zeta)$ is the Ces\`aro--Stolz Theorem, which reads as follows:
\begin{theorem}[Ces\`aro--Stolz] Let $(\Gamma_n)_{\ninn}, (\Theta_n)_{\ninn}$ be two sequences of numbers such that $\Theta_n$ is strictly increasing and tends to infinity. If the following limit exists
\be
\ds \lim_{n \to \infty} \frac{\Gamma_n - \Gamma_{n-1}}{\Theta_n - \Theta_{n-1}}
\ee then it is equal to $\lim_{n \to \infty} \Gamma_n/\Theta_n$.
\label{cstheorem}
\end{theorem}

The reader may refer to \cite{Cesaro} for the proof.

\subsection{Kooman's Theorem} \label{seckooman} Another very useful tool is an application of Kooman's Theorem to the family of $A_n(z)$'s as defined in (\ref{andef}). Kooman's Theorem, adopted for our proof, reads as follows:
\begin{theorem}[Kooman \cite{kooman1, kooman2}] Let $A$ be an $\ell \times \ell$ matrix with distinct eigenvalues. Then there exists $\epsilon > 0$ and analytic functions $U(B)$ and D(B) defined on $S_{\epsilon} = \{B: \|B - A \| < \epsilon \}$ such that \\
(1) $B = U_B D_B U_B^{-1}$, $D_B$ commutes with $A$. \\ 
(2) $U_B$ is invertible for all $B \in S_\epsilon$. \\
(3) $U_A = 1$, $D_A = A$. \\
(4) By picking a basis such that $A$ is diagonal, we can have all $D_B$ diagonal with entries being the eigenvalues of $B$.
\label{kooman}
\end{theorem}
\noindent \textit{Remark: Theorem \ref{kooman} basically follows the formulation of Theorem 12.1.7 of \cite{simon2}, except that in \cite{simon2} the statement was intended for quasi-unitary matrices. However, the same proof also holds when $A$ has distinct eigenvalues.}

The original Kooman's Theorem appeared in Theorem 1.3 of \cite{kooman1}. An application of Kooman's theorem to orthogonal polynomials was first made by Golinskii--Nevai \cite{golinskiinevai}. They applied Kooman's result to the case when $\alpha_n \to 0$ and $\sum_n \|A_{n+1} - A_n\| < \infty$ to prove that $w(\theta) > 0$ a.e. on $\T$, where $w(\theta)$ is the a.c. part of the measure.

\section{Proof of Theorem \ref{theorem1}}
The proof of Theorem \ref{theorem1} will be divided into many steps. First, we introduce a few objects and prove a lemma about them (see Lemma \ref{lemma1}). Using Lemma \ref{lemma1}, we will prove that $\lim_{n \to \infty} \Delta_n(\zeta)$ exists. Then we compute that limit explicitly and prove that the sequence $(\Delta_n(\zeta))_{\ninn}$ is of bounded variation.

\subsection{The matrix $A_n(\zeta)$ and its eigenvalues} \label{hyperbolicity} Recall the Szeg\H o recursion relations (\ref{normrec1}) and (\ref{normrec2}). Observe that they can be expressed in matrix form as follows:
\be
\bpm \vp_{n+1}(z) \\ \vp_{n+1}^*(z) \epm = (1-|\alpha_n|^2)^{-1/2} \bpm z & -\ol{\alpha_n} \\ -z \alpha_n & 1 \epm \bpm \vp_n(z) \\ \vp_n^*(z) \epm
\label{recur1}
\ee

Let
\begin{eqnarray}
A_n(z) = (1-|\alpha_n|^2)^{-1/2} \bpm z & -\ol{\alpha_n} \\ -z \alpha_n & 1 \epm
\label{andef} \\
A_\infty(z) = (1-|L|^2)^{-1/2} \bpm z & -\ol{L} \\ -z L & 1 \epm
\label{ainftydef} 
\end{eqnarray}

It is known (see Theorem 11.1.2 of \cite{simon2}) that $e^{i\theta} \in \mathbb{G}_L$ if and only if 
\begin{align}
|Tr \ainfty(e^{i \theta})|  = (1-|L|^2)^{-1/2} 2 \left|\cos \left( \ds \frac{\theta}{2} \right)\right|> 2 \end{align}

Since $\zeta$ is in the gap, $A_\infty \equiv \ainfty(\zeta)$ is hyperbolic, which implies that $\ainfty$ has two distinct eigenvalues $\lambda_1 \equiv \lambda_1(\zeta)$ and $\lambda_2 \equiv \lambda_2(\zeta)$ such that $|\lambda_1| > 1 > |\lambda_2|$ and $\lambda_2 = (\ol{\lambda_1})^{-1}$ (see Chapter 10.4 of \cite{simon2} for an introduction to the group $\mathbb{U}(1,1)$, to which $\ainfty(\zeta)$ belongs).

Let $A_n \equiv A_n(\zeta)$. Since $A_n \to \ainfty$ and $|\tr \ainfty| >2$, for some large $N_1$,
\be |\tr A_n| > 2 \quad \forall n \geq N_1\,  .
\label{N1def}
\ee
Hence, for all $n>N_1$, $A_n$ is hyperbolic and has distinct eigenvalues $\lambda_{1,n}$ and $\lambda_{2,n}$ such that $|\lambda_{1,n}|>1>|\lambda_{2,n}|$ and $\lambda_{2,n} = (\ol{\lambda_{1,n}})^{-1}$.


\subsection{$A_n(\zeta)$ and Kooman's Theorem} As seen in Section \ref{hyperbolicity} above, $A_\infty$ is hyperbolic. Hence, it has distinct eigenvalues and we can apply Kooman's Theorem (Theorem \ref{kooman}). By Kooman's Theorem, there is an open neighborhood $S_\epsilon$ around $A_\infty$ and an integer $N_2$ such that
\be A_n \in S_{\epsilon} \quad \forall n \geq N_2
\label{N2def}
\ee and there exist matrices $U_{A_n}$ and $D_{A_n}$ such that
\be
A_n = U_{A_n}\, D_{A_n}\, U_{A_n}^{-1} .
\label{un}
\ee

Perform a change of basis to make $\ainfty$ diagonal, i.e., write
\be
A_\infty = G \, D_\infty \, G^{-1}
\label{gdef}
\ee where \be
D_\infty=\bpm \lambda_1 & 0 \\ 0 & \lambda_2 \epm .
\ee By the construction of the function $D$, $D_{A_n}$ is diagonal under this new basis, so there exists a diagonal matrix 
\be
D_n = \bpm \lambda_{1,n} & 0 \\ 0 & \lambda_{2,n} \epm
\label{dndef}
\ee
such that
\be
D_{A_n} = G \, D_n \,G^{-1} .
\ee Now we define
\be G_n = U_{A_n}\, G ,
\label{Gndef}
\ee and by (\ref{un}), we have the following representation of $A_n$:
\be
A_n =  G_n \, D_n \, G_n^{-1} .
\label{eqn19b}
\ee

\subsection{The vector $w$} \label{vectorw}Let $N$ be an integer such that
\be N>\max\{N_1, N_2\} ,
\label{bigndef}
\ee where $N_1$ and $N_2$ are defined in \eqref{N1def} and \eqref{N2def} respectively. Let $w$ be the vector such that
\be
w=\bpm w_1 \\ w_2 \epm  = D_N G_N^{-1} A_{N-1} A_{N-2} \cdots A_0 \yy
\label{wdef}
\ee

We prove the following result about $w_1$ and $w_2$:
\begin{lemma}
Both $w_1$ and $w_2$ are non-zero.
\end{lemma}

\begin{proof} First of all, observe that either $w_1$ or $w_2$ must be non-zero, because both $\vp_N(\zeta)$ and $\vp_N^*(\zeta)$ are non-vanishing on $\T$, and both $D_N$ and $G_N^{-1}$ are invertible.

Now we prove $w_2 \not = 0$ by contradiction. Suppose $w_2 = 0$. Observe that $G_N w = (\vp_{N+1}(\zeta), \vp_{N+1}^*(\zeta))^T$ and $|\vp_n(\zeta)| = |\vp_n^*(\zeta)|$ on $\T$. Hence, $w_2 = 0$ implies that the matrix elements $(G_N)_{1 \, 1}$ and $(G_N)_{2 \, 1}$ satisfy
\be
|(G_N)_{1\, 1} | = |(G_N)_{2 \, 1}|
\label{eqn20a}
\ee

It will be shown later (see the discussion after (\ref{eqn41a})) that $|G_{2 \, 1} /G_{1 \, 1}| =|L| < 1$. Since $G_N \to G$, (\ref{eqn20a}) cannot be true if $N$ is sufficiently large. By a similar argument, we can also prove that $w_1 \not  =0$.

\end{proof}

\subsection{Definitions and Asymptotics of $f_{1,n}$ and $f_{2,n}$}
For $n > N$ ($N$ as defined in (\ref{bigndef})), we let
\begin{eqnarray}
P_{n} & = & \ds \prod_{k=N+1}^{n} \lambda_{1,k} \, . \label{pndef}
\end{eqnarray} Furthermore, let $f_{1,n}$ and $f_{2,n}$ be defined implicitly by the equation below:
\begin{eqnarray}
D_n G_n^{-1} G_{n-1} D_{n-1} \cdots D_{N+1} G_{N+1}^{-1} G_N w & = & P_n \bpm f_{1,n} w_1 \\ f_{2,n} w_2 \epm . \label{fndef}
\end{eqnarray} 

We are going to prove the following lemma concerning the asymptotics of $f_{1,n}$ and $f_{2,n}$:
\begin{lemma} Let $f_{1,n}$ and $f_{2,n}$ be defined as in \eqref{fndef}. The following statements hold: \\
(1) $ f_{2,n}  \to 0 {\label{f2n}}$. \\
(2) One of the following is true:
\begin{itemize}
\item (2a) There exists a constant $C$ such that $|f_{1,n}| \leq C |f_{2,n}|$. Moreover, given any $\epsilon>0$, there exist an integer $N_{\epsilon}$ and a constant $C_\epsilon$ such that
\be
|f_{2,n}| \leq C_\epsilon \left( \ds \left|\frac{\lambda_2}{\lambda_1} \right|+ \epsilon \right)^{n} , \quad \forall n \geq N_\epsilon.
\ee 
\item (2b) $|f_{2,n}/ f_{1,n}| \to 0$. Furthermore, $f_1 = \lim_{n \to \infty} f_{1,n}$ exists and it is non-zero.
\end{itemize}
\label{lemma1}
\end{lemma}

\begin{proof} \textbf{We prove statement (1) of Lemma \ref{lemma1}}.  For $n \geq N$, let the left hand side of (\ref{fndef}) be
\be
\bpm w_{1,n} \\ w_{2,n} \epm  \equiv  w(n) = D_n G_n^{-1} G_{n-1} D_{n-1} \cdots D_{N+1} G_{N+1}^{-1} G_N w \, .
\label{wndef}
\ee

First, we want to show that
\be
\|w(n+1) - D_{n+1} w(n) \| \leq C \|A_{n+1} - A_n \| |P_n| \left( |f_{1,n}| + |f_{2,n}|\right).
\label{bound3}
\ee

Note that
\be
w(n+1) - D_{n+1} w(n) = D_{n+1} \left( G_{n+1}^{-1} G_{n} - 1 \right) w(n).
\label{wnbound}
\ee

We aim to bound each of the components on the right hand side of (\ref{wnbound}). Since $U$ is analytic on $S_\epsilon$, on some compact subset of $S_\epsilon$ there exist constants $\eta_1, \eta_2>0$ such that
\be
\|G_n - G_{n-1}\|  \leq \|G\| \|U_{A_n} - U_{A_{n-1}}\|  \leq \eta_1 \| A_n - A_{n-1}\| \label{bound1} \ee and
\be
\| G_n^{-1} \| \leq  \| G^{-1}\| \| U_{A_n}^{-1}\| \leq  \eta_2 . \label{bound2}
\ee

Therefore, for $\eta = \eta_1 \eta_2$,
\be
\|G_{n+1}^{-1} G_{n} - 1\| = \| G_{n+1}^{-1} \left(  G_{n}-G_{n+1} \right) \| \leq \eta \| A_{n+1} - A_n\| .
\ee

Moreover, for $C_1=max\{|w_1|, |w_2|\}$, we have the following bounds
\begin{align}
\ds \sup_{n \geq N} \| D_{n}\| & = \sup_{n \geq N} |\lambda_{1,n}| < 2 |\lambda_1| , \\
\|w(n)\| & = \left\| \bpm f_{1,n} P_n w_1 \\ f_{2,n} P_n w_2 \epm \right\| < C_1 |P_n| \left( |f_{1,n}| + |f_{2,n}|\right) .
\end{align}

Combining all the inequalities above and applying them to (\ref{wnbound}), we have
\be
\|w(n+1) - D_{n+1} w(n) \| \leq C_2 \|A_{n+1} - A_n \| |P_n| \left( |f_{1,n}| + |f_{2,n}|\right)
\label{bound6a}
\ee where $C_2$ is a constant. This proves (\ref{bound3}). We shall see why \eqref{bound3} is useful as we prove (\ref{bound6}) and (\ref{bound7}) below. \\

Since $P_{n+1} =\lambda_{1,n+1}P_n$ and $w_{1,n}=P_n f_{1,n} w_1$ , there is a constant $C_3$ such that
\be \begin{array}{ll}
\left| f_{1,n+1} - f_{1,n} \right| & = \ds \frac{1}{|w_1|}  \left| \ds \frac{w_{1,n+1} -\lambda_{1,n+1} w_{1,n}}{P_{n+1}} \right| \\
\\
& \leq \ds \frac{1}{|w_1 P_{n+1}|} \| w(n+1) - D_{n+1} w(n)\| .
\end{array} \ee

By \eqref{bound6a}, this implies
\be \left| f_{1,n+1} - f_{1,n} \right|  \leq C_3  \| A_{n+1} - A_n\|  \left( |f_{1,n}| + |f_{2,n}|\right) .
\label{bound6}
\ee

Thus, by the triangle inequality,
\be \begin{array}{ll}
|f_{1,n+1}| & \leq |f_{1,n+1} - f_{1,n}| + |f_{1,n}| \\
& \leq \left( 1 + C_3 \|A_{n+1} - A_n\| \right) |f_{1,n}| + C_3 \|A_{n+1} - A_n\| |f_{2,n}| . 
\end{array}
\label{eqn5}
\ee

By a similar argument, one can prove that there is a constant $C_4$ such that
\be 
\left| f_{2,n+1} - \ds \frac{\lambda_{2,n}}{\lambda_{1,n}} f_{2,n}\right|  \leq C_4 \| A_{n+1} - A_n\|  \left( |f_{1,n}| + |f_{2,n}|\right) .
\label{bound7}
\ee

Similarly, by (\ref{bound7}) and the fact that $|{\lambda_{2,n}}/{\lambda_{1,n}}| < 1$,
\be
\begin{array}{ll}
|f_{2,n+1}| & \leq \left( 1 + C_4 \|A_{n+1} - A_n\| \right) |f_{2,n}| + C_4 \|A_{n+1} - A_n\| |f_{1,n}| 
\end{array}
\label{eqn6}
\ee

We add (\ref{eqn5}) to (\ref{eqn6}) to obtain
\be
|f_{1,n+1}| + |f_{2,n+1}| \leq \left(1+2C_5 \|A_{n+1} - A_n\| \right) \left( |f_{1,n}|+|f_{2,n} |\right),
\label{eqn7}
\ee where $C_5 =\max\{C_3, C_4\}$.

By applying (\ref{eqn7}) recursively, we conclude that
\be
\sup_n \left( |f_{1,n}| + |f_{2,n}| \right) < \infty .
\ee Therefore, \eqref{bound6} and \eqref{bound7} imply that $|f_{1,n+1} - f_{1,n}|$ and $|f_{2,n+1}-\lambda_{2,n} f_{2,n}/\lambda_{1,n}|$ are bounded. Furthermore, by the triangle inequality, there is a constant $C_6$ such that 
\begin{eqnarray}
|f_{1,n+1}| & \leq & |f_{1,n}| + C_6 \|A_{n+1} - A_n \|  \label{eqn8} ; \\
|f_{2,n+1}| & \leq  & \left| \ds \frac{\lambda_{2,n}}{\lambda_{1,n}} f_{2,n} \right| + C_6 \|A_{n+1} - A_n \| .
\label{eqn9}
\end{eqnarray}

By applying (\ref{eqn8}) and (\ref{eqn9}) recursively, we conclude that for any fixed $M$ such that $N \leq M \leq n$,
\begin{eqnarray}
|f_{1,n+1}| & \leq &|f_{1,M}| + C_6 \ds \sum_{j=M}^{n} \|A_{j+1} - A_j \| \, ; \\
|f_{2,n+1}| & \leq & \ds \prod_{j=M}^n \left| \frac{\lambda_{2,j}}{\lambda_{1,j}} \right| |f_{2,M}|
+C_6 \ds \sum_{j=M}^{n} \|A_{j+1} - A_j \| . \label{eqn10}
\end{eqnarray}

Without loss of generality, consider $n=2M$. Since $|\lambda_{2,n} / \lambda_{1,n}| \to |\lambda_{2}/\lambda_1| < 1$, $\prod_{j=M}^n \left| \frac{\lambda_{2,j}}{\lambda_{1,j}} \right| \to 0$ as $n \to \infty$. Moreover, $\sum_{j} \|A_{j+1}-A_j \| <\infty$ implies that $\sum_{j=M}^{n} \|A_{j+1}-A_j \| \to 0$ as $n \to \infty$.

 Therefore, $|f_{2,n}| \to 0 \text{ as } n \to \infty$. This proves (1) of Lemma \ref{lemma1}.

\medskip

\textbf{We proceed to prove statement (2) of Lemma \ref{lemma1}}.

There are two possible cases concerning $f_{1,n}$ and $f_{2,n}$: \\
Case (1): There exist a fixed integer $K$ and a constant $C$, $|f_{1,n}| \leq C |f_{2,n}|$ for all $n \geq K$. \\
Case (2): For any integer $K$ and any constant $M$, there exists an integer $n_{K,M} \geq K$ such that 
$|f_{1,n_{K,M}}| > M |f_{2,n_{K,M}}|$.

\vspace{1cm}

\textbf{Case (1)}: (\ref{bound7}) implies that for $n \geq \max\{N, K\}$, there is a constant $C_7$ such that
\be
|f_{2,n+1}| \leq \left( \left| \ds \frac{\lambda_{2,n}}{\lambda_{1,n}}\right| + C_7 \|A_{n+1} - A_n \| \right) |f_{2,n}|.
\ee

Therefore, given any $\epsilon>0$, there exist $N_{\epsilon}$ and a constant $C_\epsilon$ such that
\be
|f_{2,n}| \leq C_\epsilon \left( \ds \left|\frac{\lambda_2}{\lambda_1} \right|+ \epsilon \right)^{n} \quad \forall n \geq N_\epsilon .
\ee In other words, $f_{2,n}$ decays exponentially fast; hence, so does $f_{1,n}$. This proves (2a) of Lemma \ref{lemma1}.

\vspace{1cm}

\noindent \textbf{Case (2)}: Let $r_n = f_{2,n}/f_{1,n}$. \textbf{First, we want to show that given any $\epsilon >0$ there exists an integer $J_{\epsilon}$ such that $|r_j| < \epsilon$ for all $j \geq J_{\epsilon}$.}

First, we show that both $f_{1,n}$ and $f_{1,n+1}$ are non-zero, as \eqref{eqn16} below will involve $f_{1,n}$ and $f_{1,n+1}$ in the denominator.

By assumption, we are free to choose any $M$, so we choose an integer $M$ such that $1/M < \epsilon$. Consider any fixed pair $(K, M)$ (we will choose $K$ later in the proof). We are guaranteed the existence of an integer $n = n_{K,M} > K$ such that $|r_n| < 1/M=\epsilon$, which also implies that $f_{1,n} \not = 0$. Furthermore, by the triangle inequality and (\ref{bound6}),
\be \begin{array} {lll}
\left| \ds \frac{f_{1,n+1}}{f_{1,n}} \right| & \geq & 1- \left|\ds \frac{f_{1,n+1} - f_{1,n}}{f_{1,n}} \right| \\
\\
& \geq & 1 - C_3 \|A_{n+1} - A_n \| (1+ |r_n| ) > 0 .
\end{array}
\label{eqn12}
\ee

Thus, $f_{1,n+1}$ is also non-zero.

By the triangle inequality,
\be \begin{array}{lll}
&\left| r_{n+1} - \ds \frac{\lambda_{2,n}}{\lambda_{1,n}} r_n\right|\\
\leq & \left| \ds \frac{f_{2,n+1}}{f_{1,n+1}} - \frac{\lambda_{2,n}}{\lambda_{1,n}} \frac{f_{2,n}}{f_{1,n+1}} \right| + \ds \left| \ds \frac{\lambda_{2,n}}{\lambda_{1,n}}\right| \left| \ds \frac{f_{2,n}}{f_{1,n+1}} - \ds \frac{f_{2,n}}{f_{1,n}}\right| \\
\\
=&  \left| \ds \frac{f_{2,n+1} - (\lambda_{2,n}/\lambda_{1,n}) f_{2,n}}{f_{1,n+1}}\right| + \left| \ds \frac{\lambda_{2,n}}{\lambda_{1,n}} r_n \right| \left| \ds \frac{f_{1,n}-f_{1,n+1}}{f_{1,n+1}}  \right| .
\end{array}
\label{eqn16}
\ee

By (\ref{bound6}) and ({\ref{bound7}}), there exists a constant $C_8$ such that
\be \begin{array}{lll}
& \left| r_{n+1} - \ds \frac{\lambda_{2,n}}{\lambda_{1,n}} r_n\right| \\
\\
\leq & \ds \frac{ 1+|r_n||\lambda_{2,n}/\lambda_{1,n}|}{|f_{1,n+1}|} C_8  \| A_{n+1} - A_n\| (|f_{1,n}|+|f_{2,n}| )\\
=&  C_8 (1+|r_n||\lambda_{2,n}/\lambda_{1,n}|) \|A_{n+1} - A_n \|  \ds \frac{  |f_{1,n}|}{|f_{1,n+1}|} (1+ |r_n| ) .
\label{eqn11}
\end{array}
\ee

Furthermore, by inverting (\ref{eqn12}) one gets
\be
\left| \ds \frac{f_{1,n}}{f_{1,n+1}} \right| \leq \ds \frac{1}{1 - C_3\|A_{n+1} - A_n \| (1+ |r_n| )} \, .
\ee

Then we plug this into (\ref{eqn11}) to obtain
\be
\left| r_{n+1}\right|  \leq \left| \ds \frac{\lambda_{2,n}}{\lambda_{1,n}} r_n\right| +
 \ds \frac{C_8(1+|r_n||\lambda_{2,n}/\lambda_{1,n}|)  (1+ |r_n| )}{1 - C_3\|A_{n+1} - A_n \| (1+ |r_n| )}  \|A_{n+1} - A_n \| .
\label{eqn13}
\ee


Let $R_n$ be the second term on the right hand side of (\ref{eqn13}). Note that the quotient in front of $\|A_{n+1} - A_n \|$ is bounded. Hence, for any sufficiently large $K$, there exists $n\equiv n_{n,k} > K$ such that $|r_{n+1}|<|r_n|<\epsilon$.

Applying the same argument to $r_{n+1}$, we can prove that $|r_{n+2}| < \epsilon$. Inductively, $|r_{j}| < \epsilon$ for all large $j$. This proves $|f_{2,n}/f_{1,n}| \to 0$, the first claim of (2b) of Lemma \ref{lemma1}.\\

\textbf{It remains to show that $\lim_{n \to \infty}f_n$ exists.} We divide both sides of (\ref{bound6}) by $|f_{1,n}|$. Since $|r_n| \to 0$,
\be
\left| \ds \frac{f_{1,n+1}}{f_{1,n}} - 1 \right| \leq C  \|A_{n+1} - A_n\| \left( 1 + |r_n| \right) \to 0 .
\label{eqn18a}
\ee

Moreover, $\log$ is analytic near $1$, so in an $\epsilon$-neighborhood of $1$ there is a constant $E$ such that
\be |\log z| = |\log \zeta - \log 1| \leq E |z-1| .
\ee

By (\ref{eqn18a}),
\be
\left| \log \left(\ds \frac{f_{1,n+1}}{f_{1,n}} \right) \right| \leq C \|A_{n+1} - A_n\| .
\label{eqn19}
\ee

Therefore, the series $\sum_{j=N}^{\infty} \log \left( f_{1,j+1}/f_{1,j} \right)$ is absolutely convergent. Furthermore, as we have seen in \eqref{eqn12}, $f_{1,j} \not = 0$ for all large $j$. Thus, $\log f_{1,j}$ is finite and the following limit 
\be
\lim_{n\to \infty} \log f_{1,n+1}  = \lim_{n \to \infty}\ds \sum_{j=p}^{n} \left( \log  f_{1,j+1} - \log f_{1,j} \right) + \log f_{1,p}
\ee exists and is finite. We call the limit $\lim_{n \to \infty} f_{1,n} = f_1$. This proves the second part of (2b) and concludes the proof of Lemma {\ref{lemma1}}.

\end{proof}

\begin{proof}[\textbf{Proof of Theorem \ref{theorem1}}] By statement (2) of Lemma \ref{lemma1}, there are two possible cases:

\medskip

\noindent \textbf{First Case.} This corresponds to (2a) of Lemma \ref{lemma1}. Recall that for $n > N$,
\be
T_n \bpm 1 \\ 1 \epm = G_n P_n \bpm f_{1,n} w_1 \\ f_{2,n} w_2 \epm
\ee
and $G_n = U_{A_n} G \to G$ as $n \to \infty$. Hence, given any $\epsilon > 0$, there exists a constant $K_\epsilon$ such that
\be
\left\| T_n \bpm 1 \\ 1 \epm \right \| \leq \|G_n \| \ds \prod_{j=N}^{n} |\lambda_{1,j}| \left\| \bpm f_{1,n} w_1 \\ f_{2,n} w_2 \epm \right\| \leq K_\epsilon \left( \left|\ds \frac{\lambda_2}{\lambda_1}\right| + \epsilon \right)^n \left( |\lambda_1| +\epsilon \right)^n .
\ee This means that $|\vp_n(\zeta)|$ is exponentially decaying. As a result, $K_n(\zeta,\zeta)$ converges, $\mu(\zeta) = \lim_{n \to \infty} K_n(\zeta,\zeta)^{-1} >0$ and $\Delta_n(\zeta) \to 0$ exponentially fast. This proves claim (1) of Theorem \ref{theorem1}.

\medskip

\noindent \textbf{Second Case.} This corresponds to (2b) of Lemma \ref{lemma1}. 

\textbf{First, we compute $\lim_{n \to \infty}\Delta_n (\zeta)$ using the asymptotic expressions of $\vp_n(\zeta)$ and $\vp_n^*(\zeta)$.} By definition, $G_n \to G$. Suppose
\be
G_n = \bpm g_{1,n} & g_{1,n}' \\ g_{2,n} & g_{2,n}' \epm \to G = \bpm g_1 & g_1' \\ g_2 & g_2' \epm .
\label{gdef2}
\ee
Since $\vp_n(\zeta)$ is the first component of the vector $G_n P_n (f_{1,n} w_1, f_{2,n} w_2)^T$,
\be \begin{array}{ll}
\vp_n(\zeta) & = P_n \left( g_{1,n} f_{1,n} w_1 + g_{1,n}' f_{2,n} w_2 \right) \\
& = P_n f_{1,n} \left( g_{1,n} w_1 + g_{1,n}' r_n w_2 \right) \\
& = P_n \left( f_1 g_{1} w_1 + o(1) \right) .
\end{array}
\label{vpn}
\ee

Similarly,
\be
\vp_n^*(\zeta) = P_n \left( f_1 g_{2} w_1 + o(1) \right) .
\label{vpn*}
\ee

Since $P_n \to \infty$, both $\vp_n(\zeta)$ and $\vp_n^*(\zeta) \to \infty$. As a result, $(K_n(\zeta,\zeta))_{\ninn}$ is a positive sequence that tends to infinity. Hence, we can use the Ces\`aro--Stolz Theorem (Theorem \ref{cstheorem}). Let
\begin{align}
\Gamma_n (\zeta) & = \ol{\vp_{n+1}(\zeta)} \vp_n^*(\zeta) \label{Gammandef} \\
\Theta_n (\zeta) & = (1-\gamma) \gamma^{-1} + K_n(\zeta,\zeta) . \label{Thetandef}
\end{align}

By (\ref{vpn}) and (\ref{vpn*}),
\begin{align}
\Gamma_n(\zeta)
=  \ol{P_{n+1}} P_n \left( |f_1|^2 |w_1|^2 \ol{g_1} g_2 +o(1) \right) \label{eqn20} ; \\
\Theta_n(\zeta) -\Theta_{n-1}(\zeta)
 =   |P_n|^2 \left( |f_1|^2  |w_1|^2 |g_1|^2 + o(1) \right) .
\label{eqn21}
\end{align}

Using (\ref{eqn20}), (\ref{eqn21}) above and the fact that  $\lambda_2 = (\ol{\lambda_1})^{-1}$, we compute
\be \begin{array} {ll}
\ds \frac{\Gamma_n(\zeta) - \Gamma_{n-1}(\zeta)}{\Theta_n(\zeta) - \Theta_{n-1}(\zeta)} & = \ds \frac{\ol{P_{n+1}} P_n - \ol{P_n} P_{n-1}}{|P_n|^2} \left( \ds \frac{\ol{g_1} g_2}{|g_1|^2} + o(1) \right) \\
\\
& = \left( \ol{\lambda_{1, n+1}} - \ds \frac{1}{\lambda_{1,n}} \right) \left( \ds \frac{g_2}{g_1} + o(1) \right)\\
\\
& \to \left( \ol{\lambda_1} - \ol{\lambda_2} \right) \left( \ds \frac{g_2}{g_1}\right) .
\label{eqn23a}
\end{array}
\ee 
Since the limit in (\ref{eqn23a}) exists, $\lim_{n \to \infty} \Gamma_n(\zeta)/\Theta_n(\zeta)$ exists and is equal to the limit in (\ref{eqn23a}). It remains to compute $g_2/g_1$. Note that
\be
\bpm g_1 \\ g_2 \epm = G \bpm 1 \\ 0 \epm . \label{eqn23}
\ee By definition, $G$ is the change of basis matrix for $A_\infty$. Therefore, $g=(g_1, g_2)$ is the eigenvector of $\ainfty$ corresponding to the eigenvalue $\lambda_1$. It suffices to solve $(\ainfty - \lambda_1) g = 0$, which is equivalent to
\be
\bpm \zeta-\tau_1 & - \ol{L} \\ -\zeta L & 1 - \tau_1 \epm \bpm g_1 \\ g_2 \epm = \bpm 0 \\ 0 \epm ; \quad  \tau_1 = (1-|L|^2)^{1/2} \lambda_1 \, .
\label{eqn22} 
\ee 

Since the matrix on the left hand side of (\ref{eqn22}) has a non-zero vector in its kernel, it must have rank $1$, so the two rows are equivalent. For that reason we only have to look at the first row. Furthermore, note that we are only concerned about the ratio $g_2/g_1$, which is constant upon multiplication of $G$ by any non-zero constant; therefore, by putting $g_1 = 1$ and we deduce that
\be
\ds \frac{g_2}{g_1} = \ds \frac{\zeta - \tau_1 }{\ol{L}} .
\label{eqn41a}
\ee

Then by (\ref{eqn23a}),
\be \begin{array}{ll}
\Delta_\infty (\zeta) & = (1-|L|^2)^{1/2}  \left( \ol{\lambda_1} - \ol{\lambda_2} \right) \ds \frac{\zeta - \lambda_1(1-|L|^2)^{1/2}}{\ol{L}}
\end{array} .
\label{eqn42}
\ee

We will simplify (\ref{eqn42}) further. Let $\tau_2 = (1-|L|^2)^{1/2} \lambda_2$. Observe that $\tau_1, \tau_2$ are eigenvalues of the matrix
\be
M(\zeta) = (1-|L|^2)^{1/2} \ainfty(\zeta) = \bpm \zeta & -\ol{L} \\ -\zeta L &1 \epm .
\label{mzdef}
\ee

The characteristic polynomial of $M(\zeta)$ is
\be
f_M(y) = (\zeta-y) (1-y) - \zeta |L|^2 = y^2 - (\zeta+1)y + \zeta(1-|L|^2)
\ee and the eigenvalues of $M(\zeta)$ are
\be
y_{\pm}(\zeta)= \ds \frac{(\zeta+1) \pm \sqrt{(\zeta+1)^2 - 4\zeta(1-|L|^2)}   }{2} .
\label{ypmdef}
\ee

We do not know whether $y_+(\zeta)$ is $\tau_1$ or $\tau_2$. We decide in the following manner: 
observe that $y_{\pm}(\zeta)$ is continuous with respect to $\zeta$; hence if $|\lambda_1(\zeta_0)| >1$ for some $\zeta_0$ in the gap, we must have $|\lambda_1(\zeta)| >1$ for all $\zeta$ in the gap. Otherwise, there must be some $\zeta_1$ in the gap such that $|\lambda_1(\zeta_1)| = 1$, contradicting the hyperbolicity of $\ainfty(\zeta)$ in the gap.

Since $\zeta=1$ is in the gap, we plug it into (\ref{ypmdef}) to obtain
\be
y_{\pm}(1) = 1 \pm |L| \, .
 \ee

If we choose the branch of the square root such that $\sqrt{|L|^2}=|L|$, we have $y_+(\zeta) = \tau_1(\zeta)$ and $y_{-}(\zeta) = \tau_2(\zeta)$, and
\be \tau_1 - \tau_2 = \sqrt{(z-1)^2 + 4z|L|^2} .
\ee

Therefore,
\be
\Delta_\infty(\zeta) = \ds \overline{h(\zeta)^{1/2}} \left( \ds \frac{(\zeta-1) - h(\zeta)^{1/2}}{2\ol{L}} \right) ,
\label{deltaeq}
\ee where
\be
h(\zeta) = (\zeta-1)^2 + 4\zeta |L|^2 .
\ee This proves statement (2a) of Theorem \ref{theorem1}.

\textbf{Next, we prove statement (2b) of Theorem \ref{theorem1}.} Recall the result of Bello--L\'opez mentioned in the Introduction. Because of that, we expect $\lim_{n \to \infty} |\alpha_n(d\nu)|=|\Delta_{\infty}(\zeta)+L|=|L|$.

First, observe that for $\zeta =e^{i \theta}$,
\be
\zeta-1=\zeta^{1/2} \left( \zeta^{1/2} - \zeta^{-1/ 2} \right) = \zeta^{1/ 2} \,2 i \sin\left( \ds \frac \theta 2\right).
\ee That implies
\begin{align}
h(\zeta) & = 4 \zeta \left( |L|^2 -\sin^2\left(\frac{\theta}{2} \right) \right) , \label{h1} \\
 \ol{h(\zeta)^{1/2}} (\zeta-1)& = 4 i \sin\left(\frac{\theta}{2} \right) \sqrt{ |L|^2 -\sin^2\left(\frac{\theta}{2} \right)} \label{h2}.
\end{align}

Now we consider $\Delta_\infty(\zeta) +L$. Combining \eqref{deltaeq}, \eqref{h1} and \eqref{h2}, we have
\be 
 \Delta_\infty(\zeta) +L 
= \ds \frac{ i \, 2 \sin\left( \frac \theta 2 \right) \sqrt{|L|^2 -\sin^2\left(\frac{\theta}{2} \right) } + \left[2 \sin^2\left( \frac \theta 2 \right) - |L|^2 \right]}{\ol{L}} .
\label{eqn75}
\ee 

Since $\zeta$ is in the gap $\mathbb{G}_L$ if and only if $|L|^2 > \sin^2(\frac \theta 2)$, $\sqrt{|L|^2 - \sin^2(\theta/2)}$ is real (see Section \ref{hyperbolicity} above). Therefore, \eqref{eqn75} implies that
\begin{eqnarray}
Re \, \ol{L} \left( \Delta_\infty(\zeta) + L \right) & = & 2 \sin^2\left( \frac \theta 2 \right) - |L|^2 \label{strange1} \\
Im \, \ol{L} \left( \Delta_\infty(\zeta) + L \right) & = & 2 \sin\left( \frac \theta 2 \right) \sqrt{|L|^2 -\sin^2\left(\frac{\theta}{2} \right) }. \label{strange2}
\end{eqnarray}

Now that we have successfully separated the real and imaginary parts of $\ol{L}(\Delta_\infty(\zeta)+L)$, with a direct computation we can show that
\be
\left| \ol{L}(\Delta_\infty(\zeta) + L) \right| = |L|^2 .
\ee

It remains to compute the phase. Suppose $\ol{L} \left( \Delta_\infty(\zeta) + L \right) = |L|^2 e^{i \omega}$. $|L|^2 \cos \omega$ and $|L|^2 \sin \omega$, being the real and imaginary parts of $\ol{L} (\Delta_\infty(\zeta) + L)$ respectively, will be given by \eqref{strange1} and \eqref{strange2}. This proves statement (2b) of Theorem \ref{theorem1}.

\textbf{Now we are going to prove that $(\Delta_n(\zeta))_{\ninn}$ is of bounded variation.} 

First, we note the following estimates: \\
(1) By the definition of $A_n(\zeta)$, $\| A_n(\zeta) - A_{n-1}(\zeta) \| = O\left( |\alpha_n - \alpha_{n-1}| \right)$. \\
(2) By (\ref{bound6}), $|f_{1,n+1} - f_{1,n}| = O(\| A_{n+1}(\zeta) - A_{n}(\zeta) \|)$. \\
(3) By the definition of $G_n$ in (\ref{Gndef}), both $|g_{1,n+1} - g_{1,n}|$ and $|g_{1,n+1}' - g_{1,n}'|$ are $O(\| A_{n+1}(\zeta) - A_{n}(\zeta) \|)$. \\
(4) Since $\lambda_{1,n}$, $\lambda_{2,n}$ are the eigenvalues $A_n(\zeta)$, $|\lambda_{1,n+1} - \lambda_{1,n}|$ and $|\lambda_{2,n+1} - \lambda_{2,n}|$ are $O(|\alpha_{n+1} - \alpha_{n}|)$.\\
(5) By (\ref{eqn11}), $\left| r_{n+1} - c_n r_n \right| =O(\| A_{n+1}(\zeta) - A_{n}(\zeta) \|)$ where
\be
c_n = \ds \frac{\lambda_{2,n}}{\lambda_{1,n}} \to c = \ds \frac{\lambda_2}{\lambda_1}
\label{cndef}
\ee has norm strictly less than $1$. 
From now on, we will denote all error terms in the order of $O(|\alpha_{n}-\alpha_{n-1}|)$ as $e_n$.

\textbf{Recall that $\Delta_n(\zeta) = (1-|\alpha_n|^2)^{1/2}\Gamma_n(\zeta)/{\Theta_n(\zeta)}$. To prove that $(\Delta_n(\zeta))_{\ninn}$ is of bounded variation, we will consider $(1-|\alpha_n|^2)^{1/2}$ and $\Gamma_n(\zeta)/\Theta_n(\zeta)$ separately.}

First, note that
\be
(1-|\alpha_{n+1}|^2)^{1/2} - (1-|\alpha_{n}|^2)^{1/2} = e_{n+1} .
\ee


Recall that $f_{2,n}/f_{1,n}=r_n$. Hence, by (\ref{vpn}) and (\ref{vpn*}),
\be \begin{array}{ll}
& \ds \frac{\Gamma_n(\zeta)}{(1-\gamma)\gamma^{-1}+ K_n(\zeta,\zeta)} \\
\\
= & \ds \frac{\ol{P_{n+1}} P_n}{(1-\gamma)\gamma^{-1}+K_n(\zeta,\zeta)} \ol{f_{1,n+1}} f_{1,n} \ol{\left( g_{1,n+1} w_1 + g_{1,n+1}' r_{n+1} w_2 \right)} \left(g_{2,n} w_1 + g_{2,n}' r_n w_2 \right) \\
\\
=& \underbrace{\ds \frac{\ol{\lambda_{n+1}}|P_n|^2}{(1-\gamma)\gamma^{-1}+K_n(\zeta,\zeta)}}_{(I)} \underbrace{\ol{f_{1,n+1}} f_{1,n}}_{(II)} \underbrace{\ol{\left( g_{1,n+1} w_1 + g_{1,n+1}' r_{n+1} w_2 \right)}}_{(III)} \underbrace{\left(g_{2,n} w_1 + g_{2,n}' r_n w_2 \right)}_{(IV)} .
\end{array}
\label{1234}
\ee

\textbf{Now we will show that (I), (II), (III) and (IV) of \eqref{1234} are of bounded variation.}

We start with the easiest. For (II), note that by estimate (2) above,
\be \ol{f_{1,n+1}} f_{1,n} - \ol{f_{1,n}} f_{1,n-1} = e_n + e_{n-1} .
\ee

The next term we will estimate is (III). We start by showing that $(r_n)_{\ninn}$ is of bounded variation. Observe that
\be \begin{array}{ll}
|r_{n+1} - r_n| & \leq |c_n r_n + e_{n+1} - c_{n-1} r_{n-1} + e_n | \\
& \leq |c_n| |r_n - r_{n-1}| +  e_n + e_{n+1} \\
& \vdots \\
& \leq |c_n \dots c_1| |r_1 - r_0| + E_n + E_{n+1} ,
\end{array}
\ee where
\be
E_n = O(e_n + |c_n| e_{n-1} + |c_n c_{n-1}| |e_{n-2}| + \dots + |c_n \dots c_2| e_1) .
\label{Endef}
\ee

Hence,
\be \begin{array}{ll}
\ds \sum_{n=0}^{\infty} |r_{n+1} - r_n| & \leq |r_1 - r_0| \ds \sum_{n=1}^{\infty} |c_n \dots c_1| + \ds 2 \sum_{n=0}^{\infty} E_n
\end{array} .
\label{eqn60}
\ee

The first sum on the right hand side of (\ref{eqn60}) is finite because $|c_n| \to |c|<1$. Now we turn to the second sum. Upon rearranging,
\be
2 \sum_{n=0}^{\infty} E_n = O\left( \ds \sum_{n=0}^{\infty} e_n [1+ |c_{n+1}| + |c_{n+1} c_{n+2}| + \dots] \right) < \infty .
\label{eqn70}
\ee 

Then we observe that
\be 
\ol{\left( g_{1,n+1} w_1 + g_{1,n+1}' r_{n+1} w_2 \right)} - \ol{\left( g_{1,n} w_1 + g_{1,n}' r_{n} w_2 \right)} = e_{n+1} + O(|r_{n+1} - r_n|) .
\label{eqn62}
\ee Therefore, (III) is of bounded variation. With a similar argument we can prove that the same goes for (IV).

It remains to prove that (I) is of bounded variation. We will make use of the simple equality
\be
\ds \frac{1}{a_{n+1}} - \frac{1}{a_n} = \ds \frac{a_{n+1} - a_n}{a_{n+1} a_n} .
\ee

As a result, if $\lim_{n \to \infty}a_n = a \not = 0$ and $(a_n)_{\ninn}$ is of bounded variation, then $(1/a_n)_{\ninn}$ is also of bounded variation. \textbf{Thus, it suffices to prove that $([(1-\gamma)\gamma^{-1}+K_n(\zeta,\zeta)]/|P_n|^2)_{\ninn}$ is of bounded variation and $\lim_{n \to \infty} [(1-\gamma)\gamma^{-1}+K_n(\zeta,\zeta)]/|P_n|^2 = \mathcal{L} > 0$.}

For the convenience of computation we will define a few more objects below. First, we let
\be \Lambda_n = \begin{cases} \lambda_{1,n} & \mbox{ if } n \geq N+1 \\
1 & \mbox{ if } 0 \leq n \leq N
\end{cases} .
\ee Then by (\ref{pndef}), $P_n = \prod_{j=0}^{n} \Lambda_j$. Moreover, recall the definition of $f_{1,n}$ in (\ref{fndef}), which was only defined for $n \geq N$. For $0 \leq n \leq N$, let $f_{1,n}$ $f_{2,n}$ be defined implicitly by (\ref{vpn}) and (\ref{vpn*}). We will see later that the introduction of these objects will not affect the result of our computation.

Note that $K_n(\zeta,\zeta)$ is the summation of $n+1$ terms, so we can write
\be
\ds \frac{(1-\gamma)\gamma^{-1} + K_n(\zeta,\zeta)}{|P_n|^2} = \ds \frac{ \gamma^{-1} }{|P_n|^2} + \mathcal{T}_n \, ,
\ee where
\be
\mathcal{T}_n =  \ds \sum_{j=1}^{n} \frac{|\vp_j(\zeta)|^2}{|P_n|^2}= \ds \sum_{j=1}^{n} \ds \frac{|f_{1,j}|^2 |g_{1,j} w_1 + g_{1, j}' r_{j} w_2|^2}{|\Lambda_{j+1} \cdots \Lambda_n|^2} .
\label{mtndef}
\ee with the convention that $\Lambda_{j+1} \cdots \Lambda_n = 1$ when $j=n$.

Next, we let
\be
\mathcal{S}_n = \ds \frac{K_{n-1}(\zeta,\zeta)}{|P_{n-1}|^2} = \ds \sum_{j=0}^{n-1} \ds \frac{|f_{1,j}|^2 |g_{1,j} w_1 + g_{1, j}' r_{j} w_2|^2}{|\Lambda_{j+1} \cdots \Lambda_{n-1}|^2} .
\label{msndef}
\ee

Then
\begin{multline}
\left| \ds \frac{(1-\gamma)\gamma^{-1}+K_n(\zeta,\zeta)}{|P_n|^2} - \ds \frac{(1-\gamma)\gamma^{-1}+K_{n-1}(\zeta,\zeta)}{|P_{n-1}|^2} \right| \\ \leq \ds \frac{2(1+\gamma^{-1})}{|P_{n-1}|^2} + \left| \mT_n - \mS_n \right| .
\label{eqn61}     
\end{multline}

We will show that each of the two terms on the right hand side of (\ref{eqn61}) is summable.

Since $|\Lambda_n|^{-1} \to |\lambda_1|^{-1} < 1$,
\be
\ds \sum_{n=0}^{\infty} \ds \frac{2(1+\gamma^{-1})}{|P_n|^2} 
 = O\left( \ds \sum_{j=0}^{\infty} \ds \frac{1}{|\lambda_1|^{2j}} \right) < \infty .
\ee

Now we will go on to prove that $\mT_n - \mS_n$ is summable. Upon relabeling the indices of $\mS_n$ in \eqref{msndef}, we have
\be \begin{array}{ll}
& \mT_n - \mS_n \\
= & \ds \sum_{j=1}^{n} \left[ \ds \frac{|f_{1,j}|^2 |g_{1,j} w_1 + g_{1, j}' r_{j} w_2|^2}{|\Lambda_{j+1} \dots \Lambda_n|^2} - \ds \frac{|f_{1,j-1}|^2 |g_{1,j-1} w_1 + g_{1, j-1}' r_{j-1} w_2|^2}{|\Lambda_{j} \dots \Lambda_{n-1}|^2} \right] 
\end{array}
\label{eq60}
\ee and we will compute term by term.

Let
\be 
\epsilon_j = |g_{1,j} w_1 + g_{1,j}' r_j w_2|^2 .
\label{eqn61a}
\ee

Then by \eqref{eq60} above,
\begin{multline}
|\mT_n - \mS_n|  \leq  \underbrace{\ds \sum_{j=1}^n \ds \frac{|f_{1,j}|^2 |\epsilon_j - \epsilon_{j-1}|}{|\Lambda_{j+1} \cdots \Lambda_n|^2}}_{(I)}
+ \underbrace{\ds \sum_{j=1}^n \ds \frac{||f_{1,j}|^2 - |f_{1,j-1}|^2| \epsilon_{j-1}}{|\Lambda_{j+1} \cdots \Lambda_n|^2} }_{(II)}\\ + \underbrace{\ds \sum_{j=1}^n \ds |f_{1,j-1}|^2 \epsilon_{j-1} \left| \ds \frac{1}{|\Lambda_{j+1} \cdots \Lambda_n|^2} - \ds \frac{1}{|\Lambda_{j} \cdots \Lambda_{n-1}|^2}  \right|}_{(III)} .
\label{eq110}
\end{multline}

Now we will prove that each of the sums on the right hand side of \eqref{eq110} is summable. We will start with (II).

Recall that $|f_{1,j}-f_{1,j-1}|=O(\|A_j - A_{j-1}\|)$ and that $f_{1,j} \to f_1$. Therefore, for some constant $C$,
\begin{multline}
\ds \sum_{n=1}^{\infty} \sum_{j=1}^n \ds \frac{||f_{1,j}|^2 - |f_{1,j-1}|^2| \epsilon_{j-1}}{|\Lambda_{j+1} \cdots \Lambda_n|^2} \\< C \left( \ds \sum_{n=1}^{\infty} |f_{1,n} - f_{1,n-1}| \right) \left( \ds \sum_{j=1}^{\infty} \frac{1}{\lambda_1^{2j}}\right) < \infty .
\label{eq111}
\end{multline}

Since $g_{1,j}$, $g_{1,j}'$ and $r_j$ are all of bounded variation and their limits exist when $j$ goes to infinity, $\epsilon_j$ is of bounded variation. Hence, there exists a constant $C$ such that
\be
\ds \sum_{n=1}^{\infty} \ds \sum_{j=1}^n \ds \frac{|f_{1,j}|^2 |\epsilon_j - \epsilon_{j-1}|}{|\Lambda_{j+1} \cdots \Lambda_n|^2} <  C \left( \ds \sum_{j=1}^{\infty} |\epsilon_j - \epsilon_{j-1}| \right) \left( \ds \sum_{j=1}^{\infty} \frac{1}{\lambda_1^{2j}}\right) < \infty .
\ee

Finally, we will consider (III). Observe that
\be
\left| \ds \frac{1}{|\Lambda_{j+1} \cdots \Lambda_n|^2} - \ds \frac{1}{|\Lambda_{j} \cdots \Lambda_{n-1}|^2} \right| = \ds \frac{|\Lambda_j|^2 - |\Lambda_{n}|^2}{|\Lambda_{j} \cdots \Lambda_n|^2}
\ee and that there exists a constant $C$ independent of $j, n$ such that
\be
|\Lambda_j|^2 - |\Lambda_n|^2 = \sum_{k=j}^{n-1} \left( |\Lambda_k|^2 - |\Lambda_{k+1}|^2 \right) < C \ds \sum_{k=j}^{n-1} |\Lambda_k - \Lambda_{k+1}| .
\ee

Hence,
\be
\ds \sum_{n=1}^\infty {\ds \sum_{j=1}^n \ds \left| \ds \frac{1}{|\Lambda_{j+1} \cdots \Lambda_n|^2} - \ds \frac{1}{|\Lambda_{j} \cdots \Lambda_{n-1}|^2}  \right|} < C \ds \sum_{n=1}^{\infty}\sum_{j=1}^{n}\sum_{k=j}^{n-1} \ds \frac{|\Lambda_{k+1}- \Lambda_k|}{|\Lambda_j \cdots \Lambda_n|^2} .
\ee

Next, we count the coefficient of $|\Lambda_{k+1} - \Lambda_k|$ in the sum above. From the expression, we know that $j \leq k < n$. Therefore, the coefficient is
\begin{multline}
\ds \sum_{n=k+1}^\infty \sum_{j=1}^k \ds \frac{1}{|\Lambda_{j+1} \cdots \Lambda_n|^2} = \ds \sum_{j=1}^k \ds \sum_{n=k+1}^\infty \left( \ds \frac{1}{|\Lambda_{j+1} \cdots \Lambda_n|^2}\right)
\\ =\left( \ds \sum_{j=1}^k \ds \frac{1}{|\Lambda_2 \cdots \Lambda_{j}|^2}\right) \left( \sum_{n=k+1}^{\infty} \ds \frac{1}{|\Lambda_{k+1} \cdots \Lambda_n|^2} \right) ,
\end{multline} which is bounded above by a constant $B$ independent of $k$. This implies that (III) is summable in $n$.

As a result, $((1-\gamma)\gamma^{-1} + K_n(\zeta,\zeta)/|P_n|^2)_{\ninn}$ is of bounded variation and that implies $\lim_{n \to \infty}[(1-\gamma)\gamma^{-1}+K_n(\zeta,\zeta)]/|P_n|^2$ exists. Moreover, 
\be
\mathcal{L}=\ds \lim_{n \to \infty} \ds \frac{K_n(\zeta,\zeta)}{|P_n|^2} >\ds \lim_{n \to \infty} \ds \frac{|\vp_n(\zeta)|^2}{|P_n|^2} 
> 0 .
\ee





This concludes the proof of Theorem \ref{theorem1}.
\end{proof}

\section{proof of theorem \ref{theorem2}}

We will generalize the method developed in Theorem \ref{theorem1}. First, we define
\begin{eqnarray}
B_{k}(\zeta) & = & A(\alpha_{(k+1)p-1},z) \cdots A(\alpha_{kp},z) \label{bkdef} ;\\
B_{\infty}(\zeta) & = & A(\beta_{p-1},z) \cdots A(\beta_0,z) \label{binftydef} .
\end{eqnarray}

We need to check a few conditions concerning the $B_k(\zeta)$'s. First, note that there exists a constant $C$ such that
\be
\| B_{k+1}(\zeta) - B_k(\zeta) \| \leq C \ds \sum_{j=0}^{p-1} |\alpha_{(k+1)p+j} - \alpha_{kp+j} |
\ee

Hence,
\be \begin{array}{ll}
\ds \sum_{k=0}^{\infty} \| B_{k+1}(\zeta) - B_k(\zeta) \| & \leq C \ds \sum_{k=0}^{\infty} \sum_{j=0}^{p-1}  |\alpha_{(k+1)p+j} - \alpha_{kp+j}| \\
& = C \ds \sum_{m=0}^{\infty} |\alpha_{m+p} - \alpha_m| < \infty
\end{array}
\ee

Furthermore, since $\zeta$ is in the gap, $|\tr B_{\infty}(\zeta)| > 2$. Since $B_k(\zeta) \to B_{\infty}(\zeta)$, for all large $k$, $|\tr B_k(\zeta)| > 2$. As a result, $B_k(\zeta)$ has distinct eigenvalues $\tau_{1,k}$ and $\tau_{2,k}$ such that $|\tau_{1,k}| > 1 > |\tau_{2,k}|$ and $|\tau_{1,k}\tau_{2,k}|=1$. Moreover, $\tau_{i,k} \to \tau_{i}$, where $\tau_1, \tau_2$ are the eigenvalues of $B_{\infty}(\zeta)$.

Next, observe that for any fixed $0 \leq j \leq p-1$,
\be \begin{array}{ll}
T_{kp+j}(\zeta) & = \left( A_{kp+j}(\zeta) \cdots A_{kp}(\zeta) \right) A_{kp-1} \cdots A_0(\zeta) \\
\\
& = \left( A_{kp+j}(\zeta) \cdots A_{kp}(\zeta) \right) B_{k-1}(\zeta) B_{k-2}(\zeta) \cdots B_0(\zeta) 
\label{eqn30}
\end{array}
\ee and $A_{kp+j}(\zeta) \to A_{\infty,j}(\zeta)$, where 
\be A_{\infty, j}(\zeta) = (1-|\beta_j|^2)^{-1/2}\bpm \zeta & -\ol{\beta_j} \\ - \zeta \beta_j & 1 \epm ; 0 \leq j \leq p-1 \label{aij}
\ee

By Kooman's Theorem and a change of basis, we can express
\be
B_n(\zeta) = G_n D_n G_n^{-1}
\ee as in (\ref{eqn19b}), where $D_n$ is a diagonal matrix with entries being the eigenvalues of $B_n(\zeta)$, and $G_n \to G_\infty$, where $G_\infty$ is the matrix that diagonalizes $B_\infty(\zeta)$.

By applying an argument similar to that in Section \ref{vectorw} to the family of $B_n(\zeta)$'s, we can show that there exists a non-zero vector $w$ and an integer $N$ such that
\be
B_n(\zeta) \cdots B_0(\zeta) \bpm 1 \\ 1 \epm = G_n(\zeta) P_n \bpm f_{1,n} & 0 \\ 0 & f_{2,n} \epm \w ,
\ee where $P_n= \prod_{j=N+1}^{n} \tau_{1,j} $. Moreover, we can show that
\be
f_{1,n} \to f_1 ; \quad f_{2,n} \to f_2 ; \quad \ds \frac{f_{1,n}}{f_{2,n}} \to 0 .
\ee 

Furthermore, by (\ref{eqn30}), for each fixed $j$, we can express $T_{kp+j}(\zeta)$ as 
\be
T_{kp+j}(\zeta) v =  \left( A_{kp+j}(\zeta) \cdots A_{kp}(\zeta) \right) G_{k-1} P_{k-1}  \bpm f_{1,k-1} & 0 \\ 0 & f_{2,k-1} \epm \w
\label{eqn72}
\ee with the property that
\be
A_{kp+j}(\zeta) \cdots A_{kp}(\zeta) G_{k-1} \to A_{\infty,j}(\zeta) \cdots A_{\infty,0}(\zeta) G_\infty \equiv M_j .
\ee

Let
\be
M_j = \bpm m_{1,j} & m_{1,j'} \\ m_{2,j} & m_{2',j} \epm .
\ee

Note that for each $n$, there are two possible expressions for $T_n(\zeta) v$. We could either write it as in (\ref{eqn72}) or as follows
\be 
T_{kp+j}(\zeta)v = A_{kp+j}(\zeta) \cdots A_{(k-1)p}(\zeta) G_{k-2} P_{k-2} F_{k-2} \w
\ee The reason will be apparent later in the proof.

Consider $n=kp+j$ where $0 \leq j \leq p$. The asymptotic formulae for $\vp_n(\zeta)$ and $\vp_n^*(\zeta)$ are of the form
\begin{align}
\vp_n(\zeta) = P_{k-1} (f_1 m_{1,j} w_1 + o(1)) ; \\
\vp_n^*(\zeta) = P_{k-1} (f_1 m_{2,j+p} w_1 + o(1)) .
\end{align}

The alternate formulae for $\vp_n(\zeta)$ and $\vp_n^*(\zeta)$ are
\begin{align}
\vp_n(\zeta) = P_{k-2} ( f_1 m_{1,p+j} w_1 + o(1) ) ; \label{eqn73} \\
\vp_n^*(\zeta) = P_{k-1} (f_1 m_{2,j} w_1 + o(1)) . \label{eqn74}
\end{align}

We define $\Gamma_n(\zeta)$ and $\Theta_n(\zeta)$ as in (\ref{Gammandef}) and (\ref{Thetandef}) respectively. Then
\begin{eqnarray}
\Gamma_n(\zeta) & = & |P_{k-1}|^2 \left( |f_1|^2 |w_1|^2 \ol{m_{1,j+1}} m_{2,j} + o(1) \right) , \\
\Theta_n(\zeta) & = & |P_{k-1}|^2 \left(|f_1|^2 |w_1|^2 |m_{1,j}|^2 + o(1) \right).
\end{eqnarray}

Moreover, observe that
\be
\Gamma_{n+p}(\zeta) = |P_{k}|^2 \left( |f_1|^2 |w_1|^2 \ol{m_{1,j+1}} m_{2,j} + o(1) \right) .
\ee

Instead of $(\Gamma_{n} - \Gamma_{n-1})/(\Theta_n - \Theta_{n-1})$ in the proof of Theorem \ref{theorem1}, we compute
\be \begin{array}{ll}
& \ds \lim_{k \to \infty} \frac{\Gamma_{(k+1)p+j}(\zeta) - \Gamma_{kp+j}(\zeta)}{\Theta_{(k+1)p+j}(\zeta) - \Theta_{kp+j}(\zeta)} \\
\\
= & \ds \lim_{k \to \infty} \frac{\left( |P_k|^2 - |P_{k-1}|^2 \right)  \left( |f_1|^2 |w_1|^2 \ol{m_{1,j+1}} m_{2,j} + o(1) \right)}{ |P_{k-1}|^2 |f_1|^2 |w_1|^2 \left(|m_{1,j+p}|^2 + \dots+ |m_{1,j}|^2 + o(1) \right) } \\
\\
= & \left( |\tau_{1}|^2 - 1 \right) \ds \frac{\ol{m_{1,j+1}} m_{2,j}}{|m_{1,j+p}|^2 + \dots+ |m_{1,j}|^2} .
\end{array}
\ee

Combining with the fact that $\lim_{k \to \infty} (1-|\alpha_{kp+j}|^2)^{1/2} = (1-|\beta_j|^2)^{1/2}$, we conclude that for each fixed $0 \leq j < p$, $\lim_{k \to \infty} \Delta_{kp+j} (\zeta)$ exists.

Finally, by an argument similar to that in the proof of Theorem \ref{theorem1}, one could prove that for each fixed $j$, $(\Delta_{kp+j}(\zeta))_k$ is of bounded variation.

\section{Proof of Theorem \ref{theorem3}}
In this section, $\zeta^n \alpha_n \to L$ and $\mu(\zeta)=0$ are the only assumptions that we need. No bounded variation of the Verblunsky coefficients  is required.

Let
\begin{align}
P_n(\zeta) & = (1-|\alpha_n|^2)^{1/2} \ol{\vp_{n+1}(\zeta)} \vp_n^*(\zeta)
\end{align} and $\Theta_n(z)$ be defined as in \eqref{Thetandef}.

Note that $P_n(\zeta)/\Theta_n(\zeta) = \Delta_n(\zeta)$. Moreover, since $\mu(\zeta)=0$, $K_n(\zeta, \zeta) \to \infty$, which allows us to use the Ces\`aro--Stolz Theorem.

Let $\rho_n=(1-|\alpha_n|^2)^{1/2}$.  Since $\zeta \in \T$, we can rewrite $P_n(\zeta), P_{n-1}(\zeta)$ as follows:
\begin{align}
P_n(\zeta)  & = \rho_n \zeta^{-1} \vp_{n+1}^*(\zeta) \ol{\vp_n(\zeta)} ,  \\
P_{n-1}(\zeta) & = \rho_{n-1} \ol{\vp_{n}(\zeta)} \vp_{n-1}^*(\zeta).
\end{align}

Moreover,
\be \Theta_n(\zeta) - \Theta_{n-1}(\zeta)  = |\vp_n(\zeta)|^2
\ee and $\vp_n \not = 0$ on $\T$; therefore we could cancel $\ol{\vp_n(\zeta)}$ and obtain
\be
\ds \frac{\zeta^{n} P_n(\zeta) -\zeta^{n-1} P_{n-1}(\zeta)}{\Theta_n (\zeta)- \Theta_{n-1}(\zeta)}  = \ds \frac{\zeta^{n-1}(\rho_n \vp_{n+1}^*(\zeta) - \rho_{n-1} \vp_{n-1}^*(\zeta))}{\vp_n(\zeta)}  .
\label{eqn1e}
\ee 

By (1.5.24) and (1.5.43) in \cite{simon1} respectively,
\begin{align}
\rho_n \vp_{n+1}^*(\zeta) & = \vp_n^*(\zeta) - \alpha_n \zeta \vp_n(\zeta) ,\\
\rho_{n-1} \vp_{n-1}^*(\zeta) & = \vp_n^*(\zeta) + \alpha_{n-1} \vp_n(\zeta) .
\end{align}

Therefore, (\ref{eqn1e}) becomes
\be \begin{array}{ll}
\ds \frac{\zeta^n P_n (\zeta)- \zeta^{n-1} P_{n-1}(\zeta)}{\Theta_n (\zeta)- \Theta_{n-1}(\zeta)} & = \ds \frac{\zeta^{n-1} \left(\vp_n^*(\zeta)-\zeta \alpha_n \vp_n(\zeta) - \vp_n^*(\zeta) - \alpha_{n-1} \vp_n(\zeta)\right)}{\vp_n(\zeta)} \\
& = -(\zeta^n \alpha_n + \zeta^{n-1}\alpha_{n-1}) .
\end{array}
\label{eqn2e}
\ee

Since $\zeta^n \alpha_n \to L$, the limit of (\ref{eqn2e}) as $n \to \infty$ exists and is equal to $-2L$. Moreover, since $\zeta$ is not a pure point of $d\mu$, $\Theta_n(\zeta)$ is a strictly increasing sequence that tends to $+\infty$, so we can apply the Ces\`aro--Stolz theorem and conclude that $\zeta^n \Delta_n(\zeta) = \zeta^n P_n(\zeta)/ \Theta_n(\zeta) \to -2L$. This implies that
\be
\zeta^n \alpha_n(d\nu) = \zeta^n \alpha_n + \zeta^n \Delta_n(\zeta) \to -L .
\ee

\section{proof of corollary \ref{corollary1}}
\label{last}

First, note that $\alpha_n$ is real for all $n$, so by induction on (\ref{normrec1}) we have a closed form for $\vp_n(1)$:
\be
\vp_n(1) = \ds \prod_{j=0}^{n-1} \sqrt{\ds \frac{1-\alpha_j}{1+\alpha_j}} \in \mathbb{R} .
\label{vpne}
\ee  Moreover, since $\alpha_n \to L <0$, $\sqrt{\frac{1-\alpha_j}{1+\alpha_j}} > 1$ for large $j$, $\vp_n(1)$ is exponentially increasing towards $+\infty$. Thus, $\lim_{n \to \infty}K_n(1,1)= \infty$ and $\mu(1)=0$. By Theorem \ref{theorem1}, we have $\Delta_n(1) \to -2L$.

To prove Corollary \ref{corollary1}, we are going to show that
\be
\lim_{n \to \infty} \ds \frac{ (\Delta_n(1) + 2L)}{c_n} = -2 .
\label{eqn9e}
\ee

Observe that by (\ref{vpne}),
\be
(1-|\alpha_n|^2)^{1/2} \vp_{n+1}(1) 
= (1-\alpha_n) \vp_n(1).
\ee

Moreover, $K_n(1,1)$ is exponentially increasing. Therefore,
\be
\Delta_n(1)  + 2L =  \ds \frac{(1-\alpha_n) \vp_n(1)^2 + 2L\, K_n(1,1)}{K_n(1,1)}
+ E_n
\ee where $E_n$ is exponentially small. 

We shall use the Ces\`aro--Stolz theorem again to prove that the limit in (\ref{eqn9e}) exists and is finite. Let
\begin{align}
A_n& = c_n^{-1} \ds \left[ (1-\alpha_n) \vp_n(1)^2 + 2L \, K_n(1,1) \right] ;\\
B_n & = K_n(1,1).
\end{align} 

First, note that $B_n - B_{n-1} = \vp_n(1)^2$. Second, note that by (\ref{vpne}),
\be
(1-\alpha_{n-1}) \vp_{n-1}(1)^2 = (1 + \alpha_{n-1}) \vp_n(1)^2 .
\ee
Therefore,
\begin{multline}
A_n - A_{n-1}  = \left[ c_n^{-1} (1-\alpha_n) \vp_n(1)^2 - c_{n-1}^{-1} (1+ \alpha_{n-1}) \vp_n(1)^2 \right] \\  +c_n^{-1} (2L) K_n(1,1) - c_{n-1}^{-1} (2L) K_{n-1}(1,1) .
\label{eqn10e}
\end{multline}

 The first sum on the right hand side of (\ref{eqn10e}) is
\be
\left[ c_n^{-1} (1-L) - c_{n-1}^{-1} (1+L) -2 \right] \vp_n(1)^2 ,
\label{eqn11e}\ee while the second sum is
\be
2L \left[ c_n^{-1} \vp_n(1)^2 + (c_{n}^{-1} - c_{n-1}^{-1})K_{n-1}(1,1) \right] .
\label{eqn12e}
\ee

Combining (\ref{eqn11e}) and (\ref{eqn12e}), we have
\be
\ds \frac{A_n - A_{n-1}}{B_n - B_{n-1}} = \left[ (1+L)(c_n^{-1} - c_{n-1}^{-1}) -2  \right]  + 2L(c_{n}^{-1} - c_{n-1}^{-1})\ds \frac{ K_{n-1}(1,1)}{\vp_n(1)^2} .
\label{eqn15e}
\ee

Next, we are going to show that $\frac{ K_{n-1}(1,1)}{\vp_n(1)^2}$ exists. To do that, we use the Ces\`aro--Stolz Theorem again. Let
\begin{align}
C_n & = K_{n-1}(1,1) , \\
D_n & = \vp_n(1)^2 .
\end{align}

Recall that by (\ref{vpne}), $\vp_{n}(1)^2 = \frac{1-\alpha_n}{1+ \alpha_n} \vp_{n-1}(1)^2$. Hence,
\be
D_n - D_{n-1} = \left( \ds \frac{1-\alpha_n}{1+ \alpha_n} -1\right) \vp_{n-1}(1)^2 .
\ee

Since $C_n - C_{n-1} = \vp_{n-1}(1)^2$, we have
\be
\lim_{n \to \infty} \ds \frac{C_n - C_{n-1}}{D_n - D_{n-1}} = \lim_{n \to \infty} \left( \ds \frac{1-\alpha_n}{1+ \alpha_n} -1\right)^{-1} = \frac{1+L}{-2L} .
\ee

Therefore, $K_{n-1}(1,1)/ \vp_{n}(1)^2 = -(1+L)/2L$. By (\ref{eqn15e}) and the Ces\`aro--Stolz Theorem,
\be 
\ds \lim_{n \to \infty} \frac{A_n - A_{n-1}}{B_n - B_{n-1}} =- 2 = \lim_{n \to \infty} \frac{A_n}{B_n} .
\ee

As a result,
\be
\Delta_n(1) = -2L - 2 c_n + o\left( c_n \right) .
\ee 

This proves Corollary \ref{corollary1}. In particular, if $L=-1/2$ and $c_n=1/n$, we have the rate of convergence of $\Delta_n(1)$ being $O(1/n)$, which is clearly not exponential.

\section*{Appendix: Szeg\H o condition and bounded variation} \label{example} Both the Szeg\H o condition and bounded variation of recursion coefficients come up in the study of orthogonal polynomials very often. In this section, we will show that there is a very large class of measures with Verblunsky coefficients of bounded variation satisfying $\alpha_n \to L \not = 0$ yet failing the Szeg\H o condition (\ref{szegocondition}).


Let $d\gamma$ be a non-trivial measure on $\mathbb{R}$ such that for all $n$, $\int |x|^n d\gamma< \infty$. 
It is well-known that the family of orthonormal polynomials $(p_n(x))_{\ninn}$ obey the following recurrence relation
\be
x p_n(x) = a_{n+1} p_{n+1}(x) + b_{n+1} p_n(x) + a_n p_{n-1}(x)
\ee for $n \geq 0$. The reader should refer to \cite{sa, simon1} for details.

Remark: The reader should be reminded that the $a_n$'s and $b_n$'s in \cite{simon1} are different from those in \cite{sa}! In fact, $a_{n+1}$(\cite{simon1}) $= a_n$(\cite{sa}) and $b_{n+1}$(\cite{simon1})$=b_n$(\cite{sa}). In this paper, we are following the notations of \cite{simon1}.



Now we consider the measure $d\gamma$ on $\mathbb{R}$ which has recursion coefficients satisfying
\begin{eqnarray}
b_n \equiv 0 , \quad \quad a_n \nearrow 1 ,\\
\ds \sum_{n=1}^{\infty} |a_n-1|^2 = \infty .
\label{killipsimon}
\end{eqnarray}

This measure, supported on $[-2, 2]$, is purely a.c., and has no eigenvalues outside $[ 
-2, 2]$. Moreover, if we write $d\gamma(x)=f(x) dx$, $f(x)$ is symmetric. By the Killip--Simon Theorem \cite{killipsimon}, condition (\ref{killipsimon}) implies that such a measure fails the quasi-Szeg\H o condition, i.e.
\be
\ds \int_{[-2,2]} (4-x^2)^{1/2} \log f(x) dx = -\infty ,
\label{fails}
\ee which is weaker than the Szeg\H o condition
\be
\ds \int_{[-2,2]} (4-x^2)^{-1/2} \log f(x) dx = -\infty .
\label{fails1}
\ee

Now we consider $d\gamma_y$ supported on $[-y, y] \subset [-2,2]$, which is defined by scaling $d\gamma$
\be
d\gamma_y(x) = d\gamma\left( 2x y^{-1} \right), \quad 0 < y < 2.
\ee



Then the a.c. part of $d\gamma_y(x)$, supported on $[-y, y]$, is
\be
f_y(x) = f(2x y^{-1}) \chi_{[-y,y]} .
\ee

It is well-known that
\be
a_n(d\gamma_y) = \ds \ye a_n(d\gamma) ,\quad b_n(d\gamma_y)= \ds \ye b_n(d\gamma) .
\ee

Now we apply the inverse Szeg\H o map (see Chapter 13 of \cite{simon2}) to $d\gamma_y$ to form the probability measure $\mu_y$ on $\T$. Under this map, we have $d\mu_y(\theta) = w_y(\theta) \frac{d\theta}{2\pi}$ with
\be
w_y(\theta) = 2\pi |\sin(\theta)| f_y(2\cos \theta ) \chi_{[\theta_y, \pi - \theta_y] }(\theta) ,
\label{eqn110}
\ee where 
\be
\theta_y=\cos^{-1}\left(\frac y 2\right) \in \left(0, \frac{\pi}{2} \right) .
\ee 

For any $g$ measurable on $[-2,2]$,
\be
\ds \int g(x) d\gamma_y(x) = \ds \int g(2 \cos \theta) d\mu_y(\theta) .
\ee


By Corollary 13.1.8 of \cite{simon2}, $b_n(\gamma_y) \equiv 0$ if and only if $\alpha_{2n}(d\mu_y) \equiv 0$. Moreover, by Theorem 13.1.7 of \cite{simon2}, we know that
\be \begin{array}{ll}
a_{n+1}^2(d\gamma_y) & = (1- \alpha_{2n-1}(d\mu_y))(1-\alpha_{2n}(d\mu_y)^2) (1+ \alpha_{2n+1}(d\mu_y)) \\ & = (1- \alpha_{2n-1}(d\mu_y)) (1+ \alpha_{2n+1}(d\mu_y)) .
\end{array}
\label{eqn100}
\ee

Note that $w_y(\theta)$ is supported on two arcs, $[\theta_y, \pi-\theta_y]$ and $[\pi+\theta_y, 2\pi - \theta_y]$, and we can decompose $w_y(\theta)$ into
\be
w_y(\theta) = w_y(\theta)|_{[\theta_y, \pi-\theta_y]} + w_y(\theta)|_{[\pi+\theta_y, 2\pi - \theta_y]} .
\label{eqn101}
\ee

Moreover, because $\gamma_y(x)$ is symmetric, each of the two components on the right hand side of (\ref{eqn101}) is symmetric along the imaginary axis. Hence, we can view $d\mu_y$ as a two-fold copy of the probability measure 
\be d\nu_y(\theta) = m_y(\theta)\frac{d\theta}{2\pi}
\label{nuydef}
\ee defined on $\T$  with
\be
m_y(\theta) = 2 w_y\left( \ds \frac{\theta}{2} \right) \chi_{[2\theta_y, 2\pi -2 \theta_y]} 
\label{mydef}
\ee (this is also called the sieved orthogonal polynomials, see Example 1.6.14 of \cite{simon1}). Hence,
\be
\alpha_{2k-1}(d\mu_y) = \alpha_{k-1} (d\nu_y) .
\ee

In other words, the Verblunsky coefficients of $d\mu_y$ are
\be
0, \alpha_0(d\nu_y), 0, \alpha_1(d\nu_y), 0, \alpha_2(d\nu_y) \dots
\ee

Therefore, (\ref{eqn100}) becomes
\be
\ye^2 a_{n+1}^2(d\gamma) = (1-\alpha_{n-1}(d\nu_y)) (1+\alpha_{n}(d\nu_y))
\label{eqn102}
\ee for $n=0, 1, \dots$, with the convention that $\alpha_{-1}=-1$.

Now note that $d\nu_y$ is supported on the arc $[2\theta_y, 2\pi - 2\theta_y]$, so by the Bello-L\'opez result \cite{bello} (see also Theorem 9.9.1 of \cite{simon2}), for $a_y = \sin \left(\theta_y \right)$,
\begin{eqnarray}
\ds\lim_{n \to \infty} |\alpha_n(d\nu_y)| & = a_y \, , \\
\ds \lim_{n \to \infty} \ol{\alpha_{n+1}(d\nu_y)}\alpha_{n}(d\nu_y) & = a_y^2 \, .
\end{eqnarray} Since $\alpha_n \in \mathbb{R}$, $\alpha_n(d\nu_y)$ actually converges. Moreover, recall that $\theta_y \in (0, \frac{\pi}2) $ was defined such that $\cos(\theta_y)=\frac y 2$. Hence,
\be a_y = \sqrt{1-\cos^2(\theta_y)}=\sqrt{1-\left(\frac{y}{2}\right)^2} .
\label{aydef}
\ee

We rewrite (\ref{eqn102}) as follows
\be
\ds \left(\ds \frac{y}{2} \right)^2\frac{ a_{n+1}^2 (d\gamma)}{1-\alpha_{n-1}(d\nu_y)} - 1 = \alpha_n(d\nu_y) .
\label{eqn108}
\ee

When $n=0$, we have $\alpha_0 = (\frac {y}{2})\frac{a_{1}^2}{2} - 1 < 0$. Hence, by an inductive argument for (\ref{eqn108}) we can show that $\alpha_n < 0$ for all $n \geq 0$.

Next, we want to prove that $(\alpha_n(d\nu_y))_{\ninn}$ is of bounded variation if $(a_n(d\gamma))_{\ninn}$ is. From now on, we let $\alpha_n = \alpha_n(d\nu_y)$, $a_n = a_n(d\gamma)$ and $c=(y/2)^2 <1$.

By (\ref{eqn108}) above,
\be
\alpha_n - \alpha_{n-1} =   \ds \frac{c(a_{n+1}^2 - a_n^2)}{1-\alpha_{n-1}} + \ds \frac{c a_n^2 (\alpha_{n-1} - \alpha_{n-2})}{(1-\alpha_{n-1})(1-\alpha_{n-2})} .
\label{eqn104}
\ee

Therefore, by an inductive argument we conclude that $\sum_n (\alpha_n(d\nu_y) - \alpha_{n-1}(d\nu_y)) < \infty$ for any $0<y<2$. Hence to any monotonic sequence of $a_n \to 1$ and any $0<y<2$, there corresponds a family of $\alpha_n(d\nu_y)$'s of bounded variation that converge to $-a_y <0$.

Finally, we have to show that $m_y(\theta)$ fails the Szeg\H o condition (\ref{szegocondition}). Since $f(x)$ fails the quasi-Szeg\H o condition (\ref{fails}), it also fails the Szeg\H o condition (\ref{fails1}). Upon scaling, (\ref{fails1}) becomes
\be
\int_{y}^{-y} \left(\log f_y(x) \right) \ds \frac{1}{\sqrt{y^2 -x^2}} dx = -\infty .
\label{scaledszego}
\ee

Finally, by the Szeg\H o map and a change of variables, (\ref{scaledszego}) is equivalent to (\ref{szegocondition}).

\section*{Acknowledgements}
I would like to thank my advisor Professor Barry Simon for his time and advice; Dr. Marius Beceanu, Dr. Eric Ryckman and Dr. Maxim Zinchenko for very helpful discussions.

\end{document}